\documentclass[11pt]{article}
\usepackage[margin=1in]{geometry}
\usepackage[T1]{fontenc}
\usepackage{mathptmx}
\usepackage{amsmath,amssymb,amsthm}
\usepackage{graphicx,booktabs,xcolor,float}
\usepackage[numbers,sort&compress]{natbib}
\usepackage[hidelinks]{hyperref}
\numberwithin{equation}{section}
\newtheorem{theorem}{Theorem}[section]
\newtheorem{proposition}[theorem]{Proposition}
\newtheorem{lemma}[theorem]{Lemma}
\theoremstyle{definition}\newtheorem{definition}[theorem]{Definition}
\theoremstyle{remark}\newtheorem{remark}[theorem]{Remark}
\newcommand{\nbar}{\bar n}

\makeatletter
\renewcommand{\footnoterule}{%
  \kern-3\p@
  \noindent
  \hbox to 0.75\textwidth{%
    \leaders\hbox{-\kern0.6pt}\hfill}%
  \kern2.6\p@}
\makeatother

\title{\textbf{Fractional-Order Hybrid SIR Forecasting of COVID-19 in Bangladesh and the United States: Memory Effects and Out-of-Sample Evaluation}}
\author{\textbf{Md Babul Hossain}$^{1,}$\thanks{Corresponding author. Email: mdhossain2024@fau.edu}\hspace{0.45em} and \textbf{Mohammad Shah Alam}$^{2}$\\[2pt]
\small \itshape $^1$Department of Mathematics and Statistics, Florida Atlantic University, Boca Raton, FL 33431, USA\\
\small \itshape $^2$Department of Behavioral and Natural Sciences, Sul Ross State University, Alpine, TX, USA}
\date{}

\begin{document}
\maketitle

\begin{abstract}
\noindent Fractional-order epidemic models are often justified by memory, yet whether memory improves forecasts is rarely tested out of sample. We extend a hybrid susceptible--infectious--removed (SIR) framework, in which a daily transmission rate is reconstructed from reported cases and then forecast, to Caputo fractional order $\alpha\in(0,1]$. We prove well-posedness, conservation, and positivity, establish a linearized threshold at disease-free states, and derive a Gr\"unwald--Letnikov reconstruction that reduces exactly to the integer-order one at $\alpha=1$ and preserves nonnegativity under stated conditions. Our main theoretical result concerns the unknown initial infectious count: for $\alpha<1$ its influence tends to zero, but with an algebraic lower bound proportional to $t^{-\alpha}$, and after 332 daily steps it is about $4.4\times10^{8}$ times larger at $\alpha=0.9$ than at $\alpha=1$. In forecasting, the effect of memory depends on the setting. In Bangladesh, the integer order gives the lowest error for linear and constant-rate forecasts in 47 of 48 settings, and lowering $\alpha$ to $0.9$ significantly worsens constant-rate forecasts. In the United States, where reported cases totaled about 18\% of the U.S. population, a fractional order lowers the linear forecaster's error in all 12 comparisons in 2021 but raises the constant-rate forecaster's; no difference is significant after Holm adjustment. For U.S. seven-day case totals, the linear forecaster's MAPE fell from 16.07\% at $\alpha=1$ to 15.80\% at $\alpha=0.30$. For the neural forecaster, $\alpha=0.9$ lowers mean MAPE in 17 of 24 comparisons in Bangladesh but in only 7 of 24 in the United States. These results are consistent with memory damping over-responsive forecasts, although the source of this effect is not isolated. In neither country did fractional order consistently improve forecast accuracy. Against standard baselines, the SIR-based linear forecaster was more accurate than the classical SIR model (seven-day MAPE 10.09\% against 22.52\% in Bangladesh) and than a neural network without SIR structure, and was competitive with ARIMA: in the United States the fractional linear forecaster had the lowest error of all methods (7.54\% against 9.14\% for ARIMA), whereas in Bangladesh ARIMA had the lowest MAPE (8.45\%), with no significant difference from ARIMA in either country.


\end{abstract}

\section{Introduction}
Compartmental models describe an epidemic as a flow of people between susceptible, infectious, and removed states \cite{kermack1927,hethcote2000}. In their classical form they have no memory: the rate of change at time $t$ depends only on the state at $t$. Fractional-order models replace the time derivative with a nonlocal operator, so the rate of change depends on the whole past \cite{podlubny1999,diethelm2010}. This makes them attractive for epidemics in which earlier exposure, uneven waiting times, or delayed reporting shape what happens now. Ahmed, El-Sayed, and El-Saka studied the stability of fractional epidemic and predator--prey systems \cite{ahmed2007}. Diethelm fitted a fractional model to a dengue outbreak \cite{diethelm2013}. Angstmann, Henry, and McGann derived a fractional SIR model from a stochastic process with non-exponential recovery times \cite{angstmann2016}, and Vargas-De-Le\'on built Lyapunov functions for fractional epidemic systems \cite{vargas2015}. Neural methods have been combined with fractional operators as well, for example in fractional physics-informed neural networks \cite{pang2019}. Physics-informed networks with an embedded SIR model have also been used to identify and predict COVID-19 dynamics \cite{han2024} and to track time-varying transmission rates \cite{millevoi2024}, which is the quantity we reconstruct algebraically here.

During the COVID-19 pandemic, Caputo fractional compartmental models were fitted to data from many countries \cite{tuan2020,rezapour2020,xu2020}, including Bangladesh \cite{akter2023,das2022,paul2023}. Models with non-singular kernels, such as the Caputo--Fabrizio derivative, were proposed as alternatives \cite{baleanu2020,paul2023}. A common finding is that the fractional model fits reported data more closely than its integer-order counterpart \cite{tuan2020,rezapour2020,tassaddiq2024}. Fitted orders are usually just below one, though some are much lower; \cite{tassaddiq2024} reports $\alpha\approx0.68$ and notes that many epidemiological models adopt the Caputo operator without showing that it matters. More recently, physics-informed neural networks have been used to estimate the order together with the epidemiological parameters, giving $\alpha<1$ for COVID-19 data from Germany and Sweden and reported gains over the integer-order model \cite{zinihi2025}. We complement this work with a design in which the order is chosen on one period and tested on another. This matters because the fractional model contains the integer-order model and has one more parameter, so a closer in-sample fit is to be expected \cite{das2022}; on its own, it says nothing about prediction. Where forecasts were made, they were usually extrapolated from a single fitted window, for example two weeks ahead of one date \cite{xu2020}, rather than evaluated over many forecast origins.

Two issues limit much of this literature. The first is the choice of derivative. For a differentiable function the conformable derivative \cite{khalil2014} equals $t^{1-\alpha}f'(t)$, which is an ordinary derivative with rescaled time; it is local and carries no memory \cite{tarasov2018}. Only nonlocal operators such as the Caputo derivative can support a claim of memory, so we use the Caputo derivative throughout. The second issue is evaluation. As the studies above show, fractional models are usually judged by in-sample fit, which a model with an extra parameter will rarely lose. Large forecast evaluations during the pandemic show why out-of-sample assessment matters. In the US COVID-19 Forecast Hub, only about two-thirds of the evaluated models beat a naive baseline \cite{cramer2022}, and in the European COVID-19 Forecast Hub the accuracy of case forecasts fell as the horizon grew \cite{sherratt2023}, as it does in our results. The question that matters in practice is whether memory improves \emph{forecasts} when each forecast uses only the information available at its origin. That is the question we address.

We start from a hybrid framework developed for Bangladesh \cite{hossain2026}. A daily transmission rate is reconstructed algebraically from smoothed reported cases, its logarithm is forecast with neural and linear autoregressions, and the forecasts are propagated through the SIR equations. That framework is the integer-order case of the model studied here, so both can be compared within one forecasting framework. Changing $\alpha$ changes three things at once: the reconstructed compartments and transmission rate, the series on which the forecasters are trained, and the propagation of forecasts. The comparison therefore measures the combined effect of fractional reconstruction and propagation within that framework. We apply it to two countries with very different reporting: Bangladesh, where reported cases amount to under one percent of the population, and the United States, where they reach about 18 percent, so that the susceptible pool is depleted substantially.

Our contributions are as follows.
\begin{enumerate}
\item We formulate a Caputo fractional hybrid SIR model and prove well-posedness, conservation, positivity, and a linearized threshold at disease-free states (Section~\ref{sec:model}).
\item We derive a discrete fractional reconstruction of the compartments and transmission rate that contains the integer-order reconstruction as the case $\alpha=1$, and we prove that it conserves the total population and preserves nonnegativity under stated conditions (Section~\ref{sec:recon}).
\item We prove that for $\alpha<1$ initialization errors vanish, but no faster than algebraically, whereas they decay geometrically at $\alpha=1$. This quantifies memory and has practical consequences for initializing such models (Section~\ref{sec:memory}). Algebraic relaxation is known for fractional numerical schemes in general \cite{wangzou2023}; what we add is a direct, sign-preserving argument for this particular reconstruction, and we claim no priority for the general phenomenon.
\item We compare fractional and integer order out of sample in Bangladesh and the United States, over two periods, three infectious periods, two horizons, and two targets, using tests that allow for dependence and a diagnostic of how memory changes the forecasts, and we benchmark the SIR-based forecasters against ARIMA, a neural network without epidemic structure, and the classical SIR model (Sections~\ref{sec:design}--\ref{sec:results}).
\end{enumerate}

\section{Preliminaries}\label{sec:prelim}
\begin{definition}
For $0<\alpha<1$ and an absolutely continuous function $f$ on $[0,T]$, the Caputo derivative of order $\alpha$ is
\begin{equation}
{}^{C}\!D_t^{\alpha}f(t)=\frac{1}{\Gamma(1-\alpha)}\int_0^t (t-s)^{-\alpha}f'(s)\,ds
\end{equation}
and ${}^{C}\!D_t^{1}f=f'$.
\end{definition}
The Caputo derivative is linear and annihilates constants. On a uniform grid $t_n=nh$, the Gr\"unwald--Letnikov approximation of the Caputo derivative is \cite{podlubny1999,garrappa2018}
\begin{equation}\label{eq:gl}
{}^{C}\!D_t^{\alpha}x(t_n)\approx h^{-\alpha}\sum_{j=0}^{n}w_j\bigl(x_{n-j}-x_0\bigr)
\end{equation}
with weights
\begin{equation}\label{eq:weights}
w_0=1,\qquad w_j=\Bigl(1-\frac{\alpha+1}{j}\Bigr)w_{j-1}\quad (j\ge1)
\end{equation}

\begin{lemma}\label{lem:weights}
Let $0<\alpha<1$. Then $w_1=-\alpha$, $w_j<0$ for all $j\ge1$, and $\sum_{j\ge0}w_j=0$. The tail sums satisfy
\begin{equation}\label{eq:tail}
\tau_n:=\sum_{j>n}|w_j|=1-\sum_{j=1}^{n}|w_j|=\frac{\Gamma(n+1-\alpha)}{\Gamma(1-\alpha)\,\Gamma(n+1)}
\end{equation}
so that $\tau_n\sim n^{-\alpha}/\Gamma(1-\alpha)$ as $n\to\infty$. For $\alpha=1$, $w_1=-1$, $w_j=0$ for $j\ge2$, and $\tau_n=0$ for $n\ge1$ (with $\tau_0=1$).
\end{lemma}
\begin{proof}
Since $w_j=(-1)^j\binom{\alpha}{j}$, the weights are the Taylor coefficients of $(1-z)^\alpha$. The factor $1-(\alpha+1)/j$ equals $-\alpha<0$ for $j=1$ and is positive for $j\ge2$ because $j>1+\alpha$. Hence all $w_j$ with $j\ge1$ have the sign of $w_1=-\alpha$. Evaluating $(1-z)^\alpha$ at $z=1$ gives $\sum_j w_j=0$. The partial sums $\sum_{j=0}^{n}w_j$ are the Taylor coefficients of $(1-z)^{\alpha}/(1-z)=(1-z)^{\alpha-1}$, namely $\Gamma(n+1-\alpha)/(\Gamma(1-\alpha)\Gamma(n+1))$. Since $w_0=1$ and $w_j=-|w_j|$, this partial sum equals $1-\sum_{j=1}^n|w_j|$, and because the full sum vanishes it also equals $\sum_{j>n}|w_j|$. The asymptotic form follows from $\Gamma(n+1-\alpha)/\Gamma(n+1)\sim n^{-\alpha}$.
\end{proof}

\begin{lemma}[Extremum principle \cite{alrefai2012}]\label{lem:extremum}
Let $0<\alpha<1$ and $t_0>0$, and let $f\in C[0,t_0]\cap C^1(0,t_0]$ with $f'\in L^1(0,t_0)$. If $f(t_0)\le f(t)$ for all $t\in[0,t_0]$, then
\begin{equation}\label{eq:extremum}
{}^{C}\!D_t^{\alpha}f(t_0)=-\frac{1}{\Gamma(1-\alpha)}\left[\frac{f(0)-f(t_0)}{t_0^{\alpha}}+\alpha\int_0^{t_0}\frac{f(s)-f(t_0)}{(t_0-s)^{\alpha+1}}\,ds\right]\le 0
\end{equation}
\end{lemma}
\begin{proof}
Let $g(s)=f(s)-f(t_0)\ge0$. Since $f'\in L^1(0,t_0)$ and $f$ is continuous, $g$ is absolutely continuous on $[0,t_0]$. For $0<\varepsilon<t_0$, integration by parts gives
\[
\int_0^{t_0-\varepsilon}\frac{g'(s)}{(t_0-s)^{\alpha}}\,ds=\frac{g(t_0-\varepsilon)}{\varepsilon^{\alpha}}-\frac{g(0)}{t_0^{\alpha}}-\alpha\int_0^{t_0-\varepsilon}\frac{g(s)}{(t_0-s)^{\alpha+1}}\,ds
\]
Since $f'$ is continuous at $t_0$, it is bounded near $t_0$ by some $M$, so $0\le g(t_0-\varepsilon)\le M\varepsilon$ and the first term tends to zero. In the last integral the integrand is nonnegative and bounded by $M(t_0-s)^{-\alpha}$ near $t_0$, so the integral converges as $\varepsilon\to0$. Letting $\varepsilon\to0$ and dividing by $\Gamma(1-\alpha)$ gives \eqref{eq:extremum}, and both terms in the bracket are nonnegative.
\end{proof}

\section{The fractional hybrid SIR model}\label{sec:model}
Let $S$, $I$, and $R$ denote susceptible, infectious, and removed individuals in a population of size $N$. Throughout, $R_t$ denotes the removed compartment and $R_0$ its initial value; reproduction numbers are written in calligraphic type, $\mathcal{R}_0$ for the basic and $\mathcal{R}_t$ for the model-based effective reproduction number. We consider
\begin{equation}\label{eq:model}
\begin{aligned}
{}^{C}\!D_t^{\alpha}S&=-\beta(t)\frac{SI}{N}\\
{}^{C}\!D_t^{\alpha}I&=\beta(t)\frac{SI}{N}-\gamma I\\
{}^{C}\!D_t^{\alpha}R&=\gamma I
\end{aligned}
\end{equation}
where $0<\alpha\le1$, $\gamma>0$, and $\beta\in C^1[0,T]$ with $0\le\beta(t)\le\beta_{\max}$. We measure time in days and treat it as dimensionless, $t/(1\ \text{day})$, so that $\beta$ and $\gamma$ are dimensionless model parameters and the daily grid has step $h=1$. For $\alpha=1$, the choice $\gamma=1/D$ makes $D$ the mean infectious period in days. For $\alpha<1$ this interpretation is lost: the residence time in the infectious class implied by ${}^{C}\!D_t^{\alpha}I=-\gamma I$ has survival function $E_\alpha(-\gamma t^\alpha)$, which decays like $t^{-\alpha}$ and therefore has infinite mean (cf.\ \cite{angstmann2016}). We use $\gamma=1/D$ as the baseline and treat $D\in\{7,10,14\}$ as a label for the removal scenario, not as a mean duration; Section~\ref{sec:sens} repeats the analysis with the scaling $\gamma=(1/D)^\alpha$, which is common in fractional epidemic models \cite{diethelm2013}.

\begin{theorem}[Well-posedness, conservation, positivity]\label{thm:wp}
Let $S_0,I_0,R_0\ge0$ with $S_0+I_0+R_0=N$. For every $T>0$, system \eqref{eq:model} has a unique continuous solution on $[0,T]$. The solution satisfies $S+I+R=N$ and $S,I,R\ge0$ on $[0,T]$.
\end{theorem}
\begin{proof}
Write $x=(S,I,R)$ and system \eqref{eq:model} as ${}^{C}\!D_t^{\alpha}x=F(t,x)$. Because $\beta\in C^1[0,T]$, the function $F$ is continuously differentiable in $(t,x)$ and in particular locally Lipschitz in $x$. By the existence, uniqueness, and continuous-dependence theory for Caputo systems \cite[Ch.~6]{diethelm2010} (see also \cite{li2018}), a unique continuous solution exists on a maximal interval $[0,t^\dagger)$ and depends continuously on the initial value and on perturbations of $F$. For a continuously differentiable right-hand side, the regularity theory in the same chapter shows that the solution is continuously differentiable on $(0,t^\dagger)$, with a derivative that may behave like $t^{\alpha-1}$ near $t=0$ and is therefore integrable. These are the properties required by Lemma~\ref{lem:extremum}, and they hold equally for the perturbed systems used below, whose right-hand sides differ from $F$ by a constant.

\emph{Conservation.} Adding the three equations and using linearity gives ${}^{C}\!D_t^{\alpha}(S+I+R)=0$, so $S+I+R\equiv N$.

\emph{Quasi-positivity.} For every $t$ and every $x\ge0$, the $i$-th component of $F$ is nonnegative whenever $x_i=0$: on $S=0$ the first component is $0$, on $I=0$ the second is $0$, and on $R=0$ the third is $\gamma I\ge0$.

\emph{Strict positivity of a perturbed system.} For $\varepsilon\in(0,1]$ let $x^\varepsilon$ solve
\[
{}^{C}\!D_t^{\alpha}x^\varepsilon=F(t,x^\varepsilon)+\varepsilon\mathbf{1},\qquad x^\varepsilon(0)=x(0)+\varepsilon\mathbf{1}
\]
where $\mathbf{1}=(1,1,1)$. We claim that $x^\varepsilon(t)>0$ componentwise on its interval of existence. Since $x^\varepsilon(0)>0$ and $x^\varepsilon$ is continuous, otherwise $t^*=\inf\{t>0:\ x_i^\varepsilon(t)\le0\text{ for some }i\}$ is positive, and there is an index $i$ with $x_i^\varepsilon(t^*)=0$, $x^\varepsilon(t^*)\ge0$, and $x_i^\varepsilon(t)>0$ for $t<t^*$. Thus $x_i^\varepsilon$ attains its minimum over $[0,t^*]$ at $t^*$, and Lemma~\ref{lem:extremum} gives ${}^{C}\!D_t^{\alpha}x_i^\varepsilon(t^*)\le0$; for $\alpha=1$, where Lemma~\ref{lem:extremum} does not apply, the same conclusion follows from the ordinary derivative, since a differentiable function attaining its minimum over $[0,t^*]$ at $t^*$ has $(x_i^\varepsilon)'(t^*)\le0$. By quasi-positivity, however, ${}^{C}\!D_t^{\alpha}x_i^\varepsilon(t^*)=F_i(t^*,x^\varepsilon(t^*))+\varepsilon\ge\varepsilon>0$, a contradiction.

\emph{A priori bound.} The sum $\sigma^\varepsilon=x_1^\varepsilon+x_2^\varepsilon+x_3^\varepsilon$ satisfies ${}^{C}\!D_t^{\alpha}\sigma^\varepsilon=3\varepsilon$, so $\sigma^\varepsilon(t)=N+3\varepsilon+3\varepsilon t^\alpha/\Gamma(\alpha+1)$. Together with positivity, each component of $x^\varepsilon$ is bounded on bounded intervals, so $x^\varepsilon$ exists on $[0,T]$ and remains in the compact set $K=\{y\ge0:\ y_1+y_2+y_3\le N+3+3T^\alpha/\Gamma(\alpha+1)\}$ for all $\varepsilon\in(0,1]$.

\emph{Limit $\varepsilon\to0$.} Since $F$ is Lipschitz on $K$, continuous dependence gives $x^\varepsilon\to x$ uniformly on every compact subinterval of $[0,t^\dagger)\cap[0,T]$. Hence $x\ge0$ there. Combined with $S+I+R=N$, the solution stays in the compact set $\{y\ge0:\ y_1+y_2+y_3=N\}$, so it cannot blow up and $t^\dagger>T$. Thus the solution exists, is unique, conserves $N$, and is nonnegative on $[0,T]$.
\end{proof}

\begin{proposition}[Linearized threshold at disease-free states]\label{thm:stab}
Let $\beta>0$ be constant. Every point $(S_*,0,N-S_*)$ with $0\le S_*\le N$ is an equilibrium of \eqref{eq:model}. Linearizing at such a point gives, for the infectious component,
\begin{equation}\label{eq:linI}
{}^{C}\!D_t^{\alpha}y=\lambda y,\qquad \lambda=\frac{\beta S_*}{N}-\gamma
\end{equation}
whose solution $y(0)E_\alpha(\lambda t^\alpha)$ tends to zero if $\beta S_*/(\gamma N)<1$ and grows without bound if $\beta S_*/(\gamma N)>1$. For fixed $\beta$ and $\gamma$, the threshold $\beta S_*/(\gamma N)=1$ does not depend on $\alpha$. Because the disease-free equilibria form a continuum, none of them is asymptotically stable. The system has no endemic equilibrium.
\end{proposition}
\begin{proof}
At an equilibrium all right-hand sides of \eqref{eq:model} vanish; the third equation gives $0=\gamma I_*$, and since $\gamma>0$, $I_*=0$, so there is no equilibrium with $I_*>0$. Conversely, if $I=0$, all right-hand sides vanish, so every point $(S_*,0,N-S_*)$ is an equilibrium. At such a point the Jacobian of the right-hand side has the single nonzero column $\partial F/\partial I=(-\beta S_*/N,\ \lambda,\ \gamma)^{\top}$, so its eigenvalues are $0$, $0$, and $\lambda$; the zero eigenvalues correspond to directions along the continuum of equilibria. The scalar equation \eqref{eq:linI} has the solution $y(0)E_\alpha(\lambda t^\alpha)$ \cite{diethelm2010}. For $\lambda<0$, $E_\alpha(\lambda t^\alpha)\to0$, like $t^{-\alpha}/(|\lambda|\Gamma(1-\alpha))$ when $\alpha<1$; for $\lambda>0$, it grows exponentially. Finally, every neighbourhood of a disease-free equilibrium contains other disease-free equilibria, and solutions starting there remain there, so no disease-free equilibrium attracts all nearby solutions.
\end{proof}

\begin{remark}\label{rem:ml}
For fixed $\beta$ and $\gamma$, the fractional order does not move the linearized threshold, which is the effective reproduction number $\beta S_*/(\gamma N)=1$. Under the alternative scaling $\gamma=(1/D)^\alpha$ of Section~\ref{sec:sens}, however, $\gamma$ itself varies with $\alpha$, and so does the threshold ratio. For fixed $\beta$ and $\gamma$ below the invasion threshold, infectious perturbations decay algebraically, like $t^{-\alpha}$, for $0<\alpha<1$ and exponentially for $\alpha=1$.

For a fully susceptible population, $S_*=N$, the threshold ratio is $\mathcal{R}_0=\beta/\gamma$. For the integer-order model this is the next-generation value: with the single infected compartment, the new-infection and transition terms are $\mathcal{F}=\beta SI/N$ and $\mathcal{V}=\gamma I$, so at $S=N$ the next-generation matrices are $F=[\beta]$ and $V=[\gamma]$, and $\mathcal{R}_0=\rho(FV^{-1})=\beta/\gamma$. For $\alpha<1$ we use the same ratio only as the linearized invasion threshold of Proposition~\ref{thm:stab}, since $\lambda=\gamma(\mathcal{R}_0-1)$ at $S_*=N$ has the sign of $\mathcal{R}_0-1$; we do not attach to it the generation-counting interpretation of the integer-order case. Holding the other parameter fixed, the normalized (local) sensitivity of $\mathcal{R}_0$ is $(\partial\mathcal{R}_0/\partial\beta)(\beta/\mathcal{R}_0)=+1$ with respect to $\beta$ and $(\partial\mathcal{R}_0/\partial\gamma)(\gamma/\mathcal{R}_0)=-1$ with respect to $\gamma$; to first order, a small relative change in either parameter changes $\mathcal{R}_0$ by the same relative amount, with opposite signs.
\end{remark}

\section{Discrete fractional reconstruction}\label{sec:recon}
Let $n_t$ denote daily reported cases and $\nbar_t$ their causal seven-day trailing average, which uses only reports available by day $t$ \cite{hossain2026}. As in the integer-order framework, we identify model incidence with the smoothed reported count:
\begin{equation}\label{eq:incid}
\beta_t\frac{S_tI_t}{N}=\nbar_{t+1}
\end{equation}
Applying \eqref{eq:gl} with $h=1$ to \eqref{eq:model}, evaluating the right-hand side at $t$, and using \eqref{eq:incid} gives the explicit recursions
\begin{equation}\label{eq:recon}
\begin{aligned}
S_{t+1}&=S_0-\nbar_{t+1}-\sum_{j=1}^{t+1}w_j\bigl(S_{t+1-j}-S_0\bigr)\\
I_{t+1}&=I_0+\nbar_{t+1}-\gamma I_t-\sum_{j=1}^{t+1}w_j\bigl(I_{t+1-j}-I_0\bigr)\\
R_{t+1}&=R_0+\gamma I_t-\sum_{j=1}^{t+1}w_j\bigl(R_{t+1-j}-R_0\bigr)
\end{aligned}
\end{equation}
initialized with $I_0=\nbar_0$, $S_0=N-\nbar_0$, and $R_0=0$. Whenever $S_t,I_t>0$, the reconstructed transmission rate is
\begin{equation}\label{eq:beta}
\beta_t=\frac{N\,\nbar_{t+1}}{S_tI_t}
\end{equation}
Every quantity at time $t$ depends only on $\nbar_0,\dots,\nbar_{t+1}$, so the reconstruction is causal. In particular, $\beta_t$ can be computed only once the report for day $t+1$ is available, so at a forecast origin $T$ the latest available rate is $\beta_{T-1}$ (Section~\ref{sec:design}). In the fractional setting, \eqref{eq:incid} should be read as a reported-case forcing approximation: the term $\beta_tS_tI_t/N$, which in the integer-order model is the number of new infections per day, is identified with smoothed daily reports, and the fractional recursion then determines how these reports accumulate into prevalence.

\begin{proposition}[Consistency with integer order]\label{prop:alpha1}
For $\alpha=1$, the recursions \eqref{eq:recon} reduce to $S_{t+1}=S_t-\nbar_{t+1}$, $I_{t+1}=(1-\gamma)I_t+\nbar_{t+1}$, and $R_{t+1}=R_t+\gamma I_t$, the integer-order reconstruction of \cite{hossain2026}.
\end{proposition}
\begin{proof}
For $\alpha=1$, Lemma~\ref{lem:weights} gives $w_1=-1$ and $w_j=0$ for $j\ge2$, so each sum in \eqref{eq:recon} equals $-(x_t-x_0)$ and the recursions collapse as stated.
\end{proof}

\begin{proposition}[Discrete conservation and positivity]\label{prop:discrete}
The reconstruction satisfies $S_t+I_t+R_t=N$ for all $t$. If $0<\gamma\le\alpha\le1$ and $\nbar_t\ge0$, then $I_t\ge0$ and $R_t\ge0$ for all $t$, and
\begin{equation}\label{eq:Ipos}
I_{t+1}=\nbar_{t+1}+(\alpha-\gamma)I_t+\sum_{j=2}^{t+1}|w_j|\,I_{t+1-j}+\tau_{t+1}I_0
\end{equation}
Moreover, $S_{t+1}\ge0$ whenever $\nbar_{t+1}\le S_0\tau_{t+1}+\sum_{j=1}^{t+1}|w_j|S_{t+1-j}$.
\end{proposition}
\begin{proof}
Let $N_t=S_t+I_t+R_t$. Summing \eqref{eq:recon} gives $N_{t+1}-N_0=-\sum_{j=1}^{t+1}w_j(N_{t+1-j}-N_0)$, which vanishes by induction on $t$. For the infectious recursion, write $w_j=-|w_j|$ and use $\sum_{j=1}^{t+1}|w_j|=1-\tau_{t+1}$ from Lemma~\ref{lem:weights}:
\[
I_{t+1}=\nbar_{t+1}-\gamma I_t+\sum_{j=1}^{t+1}|w_j|I_{t+1-j}+\tau_{t+1}I_0
\]
Separating the term $|w_1|I_t=\alpha I_t$ gives \eqref{eq:Ipos}. When $\gamma\le\alpha$ every coefficient is nonnegative, so $I_t\ge0$ by induction. The same rearrangement gives $R_{t+1}=\gamma I_t+\sum_{j=1}^{t+1}|w_j|R_{t+1-j}+\tau_{t+1}R_0\ge0$ and $S_{t+1}=S_0\tau_{t+1}+\sum_{j=1}^{t+1}|w_j|S_{t+1-j}-\nbar_{t+1}$.
\end{proof}

In the application $\gamma\le1/7<0.70\le\alpha$, and under the alternative scaling of Section~\ref{sec:sens} $\gamma=(1/D)^\alpha\le(1/7)^{0.7}<0.26$, so the conditions for $I_t\ge0$ and $R_t\ge0$ hold for every order considered. The condition for $S_t$ is an incidence budget, which we checked numerically: $S_t/N\ge0.9906$ in every reconstruction, and in all 3{,}900 simulated forecast trajectories (every origin, order, infectious period, forecaster, and scaling) $S_t/N\ge0.9907$ and $I_t>0$, with no failures. On the Bangladesh record the reconstruction satisfies $S_t+I_t+R_t=N$ to within $5\times10^{-7}$ persons and $\beta_t>0$ throughout, and for $\alpha=1$ it agrees with the integer-order reconstruction of \cite{hossain2026} to a relative difference of $3\times10^{-15}$.

\section{Memory: persistence of initialization effects}\label{sec:memory}
The initial infectious count $I_0$ is not observed; it is set to $\nbar_0$. In the integer-order model an error $\delta$ in $I_0$ decays geometrically, $I_t^{(\delta)}-I_t=(1-\gamma)^t\delta$ \cite{hossain2026}, so a long initialization period makes it negligible. The fractional model behaves differently.

\begin{theorem}[Slow forgetting of initialization errors]\label{thm:memory}
Let $0<\gamma\le\alpha<1$ and replace $I_0$ by $I_0+\delta$ in \eqref{eq:recon}, keeping the incidence series fixed. Then $e_t:=I_t^{(\delta)}-I_t$ has the sign of $\delta$ for all $t$, and:
\begin{enumerate}
\item[(i)] (algebraic lower bound) for $t\ge1$,
\begin{equation}\label{eq:memory}
|e_t|\ \ge\ |\delta|\,\tau_t=|\delta|\,\frac{\Gamma(t+1-\alpha)}{\Gamma(1-\alpha)\,\Gamma(t+1)}\ \sim\ \frac{|\delta|\,t^{-\alpha}}{\Gamma(1-\alpha)}
\end{equation}
\item[(ii)] (convergence) $|e_t|\le(1-\gamma)|\delta|$ for $t\ge1$, and $e_t\to0$ as $t\to\infty$.
\end{enumerate}
\end{theorem}
\begin{proof}
The recursion \eqref{eq:Ipos} is linear in $(I_0,\dots,I_t)$ for a fixed incidence series. Subtracting its two versions gives
\[
e_{t+1}=(\alpha-\gamma)e_t+\sum_{j=2}^{t+1}|w_j|\,e_{t+1-j}+\tau_{t+1}\,\delta,\qquad e_0=\delta
\]
All coefficients are nonnegative because $\gamma\le\alpha$, so by induction every $e_t$ has the sign of $\delta$. We may therefore take $\delta>0$.

(i) All terms on the right are nonnegative; dropping all but the last gives $e_{t+1}\ge\tau_{t+1}\delta$. The closed form and asymptotics of $\tau_t$ are given by Lemma~\ref{lem:weights}.

(ii) By Lemma~\ref{lem:weights}, $\sum_{j=2}^{t+1}|w_j|=1-\alpha-\tau_{t+1}$. If $0\le e_s\le\delta$ for all $s\le t$, then
\[
e_{t+1}\le\bigl[(\alpha-\gamma)+(1-\alpha-\tau_{t+1})\bigr]\delta+\tau_{t+1}\delta=(1-\gamma)\delta
\]
so by induction $0\le e_t\le(1-\gamma)\delta\le\delta$ for all $t\ge1$. Let $L=\limsup_{t\to\infty}e_t\in[0,\delta]$ and fix $\varepsilon>0$. Choose $J$ with $\tau_J<\varepsilon$ and then $t_0$ such that $e_s\le L+\varepsilon$ for $s\ge t_0$. For $t\ge t_0+J$, split the sum at $J$:
\[
e_{t+1}\le(\alpha-\gamma)(L+\varepsilon)+\sum_{j=2}^{J}|w_j|(L+\varepsilon)+\tau_J\,\delta+\tau_{t+1}\delta\le(1-\gamma)(L+\varepsilon)+\varepsilon\delta+\tau_{t+1}\delta
\]
Taking $\limsup$ and using $\tau_{t+1}\to0$ gives $L\le(1-\gamma)(L+\varepsilon)+\varepsilon\delta$. Letting $\varepsilon\to0$ yields $\gamma L\le0$, so $L=0$ and $e_t\to0$.
\end{proof}

\begin{figure}[H]
\centering\includegraphics[width=0.9\linewidth]{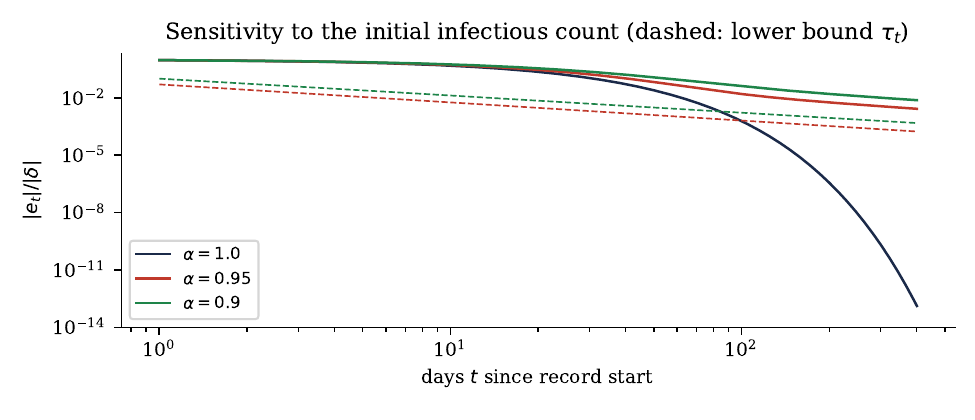}
\caption{Effect of a unit change in the initial infectious count on reconstructed prevalence, computed from the exact error recursion (solid), and the lower bound \eqref{eq:memory} (dashed). The integer-order effect decays geometrically; the fractional effects decay much more slowly.}\label{fig:memory}
\end{figure}

\begin{remark}
Part (i) is a lower bound only: it shows that initialization errors cannot vanish faster than $t^{-\alpha}$, but it does not establish the exact rate. A matching upper bound is not proved here, although the numerical decay in Figure~\ref{fig:memory} follows the slope of the bound; exact polynomial decay rates for related fractional numerical schemes are established in \cite{wangzou2023}. For $\alpha=1$ the forcing term vanishes ($\tau_t=0$) and the exact geometric decay $e_t=(1-\gamma)^t\delta$ is recovered.
\end{remark}

The error recursion involves only $\alpha$, $\gamma$, and $t$, not the incidence series, so the values below are properties of the scheme rather than estimates from Bangladesh data. Figure~\ref{fig:memory} shows the normalized error $|e_t|/|\delta|$ for $\gamma=1/14$. After 332 steps, the length of the initialization period before model fitting in \cite{hossain2026}, the effect, computed from the exact error recursion in the proof of Theorem~\ref{thm:memory}, is $2.1\times10^{-11}$ at $\alpha=1$, equal to $(13/14)^{332}$, but $3.2\times10^{-3}$ at $\alpha=0.95$ and $9.1\times10^{-3}$ at $\alpha=0.9$, about $4.4\times10^{8}$ times the integer-order value ($9.1\times10^{-3}$ against $2.1\times10^{-11}$) and above the bounds $\tau_{332}=2.1\times10^{-4}$ and $5.7\times10^{-4}$. For fractional models, a long initialization period does not by itself make the unknown initial state irrelevant, and the initial state should be treated as an uncertain input.

\begin{figure}[H]
\centering\includegraphics[width=0.9\linewidth]{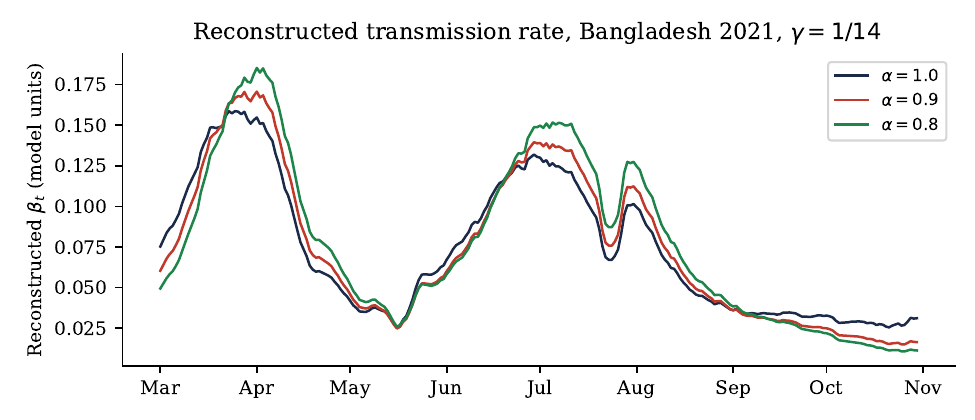}
\caption{Transmission rates reconstructed by \eqref{eq:recon}--\eqref{eq:beta} from Bangladesh data for $\alpha=1,0.9,0.8$, with $D=14$. The curve for $\alpha=1$ is the integer-order reconstruction.}\label{fig:beta}
\end{figure}

\section{Forecasting design}\label{sec:design}
\subsection{Data}
We use the compiled daily record of laboratory-confirmed COVID-19 cases in Bangladesh from 3 April 2020 to 8 January 2022 (646 days), the population $N=169{,}828{,}911$, and the causal smoothing of \cite{hossain2026}, $\nbar_t=\min(t+1,7)^{-1}\sum_{k=0}^{\min(t,6)}n_{t-k}$, which uses only reports available by day $t$. Every daily count in the record is positive (the minimum is 5, and no day has zero reports). All targets are therefore positive, MAPE is well defined, and no zero-count rule is needed. The data, their check against official monthly totals, and their properties are described in \cite{hossain2026}. Reported cases understate infections \cite{bhuiyan2022}, so all compartments should be read as effective reported-case quantities.

For the United States we use the daily national case series compiled by The New York Times \cite{nyt2021} over the same 646 days, with the 2020 Census resident population $N=331{,}449{,}281$ \cite{census2021}. Reported cases in this window sum to 59.6 million, about 18\% of the population, against under 1\% in Bangladesh. The series has a strong weekday pattern: mean daily reports are 54{,}943 on Sundays and 111{,}682 on Mondays. The cumulative series contains one downward correction inside the window ($-23{,}999$ cases on 4 June 2021). We replace it by the mean of the three preceding and three following days (15{,}825), which changes the total by less than 0.1\%; the smallest daily count is then 2{,}922, so all targets remain positive. The corrected value enters training only at origins from 10 June 2021 onward, when the six neighbouring days are already observed. Through the seven-day smoothing it also enters the evaluation targets of origins whose forecast windows cover 4 to 10 June, namely the 3 June origin at both horizons and the 27 May origin at fourteen days. Excluding the affected forecast origins preserved the direction of all 12 MAPE differences for each of the linear and constant-rate forecasters, and none reached statistical significance after Holm adjustment (Table~\ref{tab:usexcl} in Appendix~A).

\subsection{Forecasters}
At origin $T$ the reconstructed states $S_T$, $I_T$, and $R_T$ are known, but the rate is known only up to $\beta_{T-1}$, because \eqref{eq:beta} requires the incidence of day $T+1$; the first forecast is for day $T+1$. Let $b_t=\log\beta_t$. As in \cite{hossain2026}, each method forecasts $b_T,\dots,b_{T+H-1}$ recursively from seven lags: (i) a ridge-regularized linear autoregression (penalty one, standardized inputs); (ii) a constant rate equal to the mean of the last seven reconstructed rates; and (iii) an ensemble of ten feed-forward networks with two hidden layers of eight and four tanh units, fitted by L-BFGS with penalty $0.01$, at most 20{,}000 iterations and $10^6$ function evaluations, tolerance $10^{-7}$, and seeds 0--9 (scikit-learn 1.8.0, Python 3.12); the pointwise median rate of the ensemble is used. Training uses rates from 1 March 2021 (1 May 2020 for the 2020 period) up to $T-1$, with all lags inside the window. At each origin we standardize every input column and the target with the mean and standard deviation of the training pairs available at that origin only, and we transform forecasts back before the next recursive step. If an input column has zero training variance, we only centre it, setting its scale to one.

\subsection{Fractional propagation}
Starting from the reconstructed history up to $T$, each forecast rate is inserted into the discrete fractional system: for $t=T,\dots,T+H-1$,
\begin{equation}\label{eq:prop}
\begin{aligned}
\hat n_{t+1}&=\hat\beta_t\frac{\hat S_t\hat I_t}{N}\\
\hat S_{t+1}&=S_0-\hat n_{t+1}-\sum_{j=1}^{t+1}w_j\bigl(\hat S_{t+1-j}-S_0\bigr)\\
\hat I_{t+1}&=I_0+\hat n_{t+1}-\gamma\hat I_t-\sum_{j=1}^{t+1}w_j\bigl(\hat I_{t+1-j}-I_0\bigr)
\end{aligned}
\end{equation}
with $\hat S_s=S_s$ and $\hat I_s=I_s$ for $s\le T$. Unlike in the integer-order model, each step depends on the whole history, both reconstructed and simulated.

\subsection{Baseline models}\label{sec:baseline_models}
To place the SIR-based forecasters among standard alternatives, we also evaluate three baselines at the same origins, with the same training windows. The first is the classical SIR model, represented by the constant-rate forecaster at $\alpha=1$, which keeps the transmission rate at its last weekly mean. The second is an ARIMA$(p,1,q)$ model \cite{boxjenkins2015} for $\log\nbar_t$, with $p\le3$, $q\le2$ and an optional drift. At each origin the order is selected by the Akaike information criterion \cite{akaike1974}, in the spirit of automatic ARIMA selection \cite{hyndman2008}, all candidates are fitted by conditional sum of squares on the same effective sample, and candidates with a non-stationary autoregressive or non-invertible moving-average part are discarded. The third (ANN only) is the neural ensemble described above, trained directly on seven lags of $\log\nbar_t$ rather than of $\log\beta_t$, so that the two neural forecasters share the same network architecture and differ in what they forecast: incidence directly, rather than the transmission rate propagated through the SIR equations. Forecasts of the logarithm are transformed back by exponentiation, as for the SIR-based forecasters. The baseline comparison uses $D=14$ and the 2021 origins.

\subsection{Choice of order and evaluation}
We consider $\alpha\in\{0.70,0.75,\dots,1.00\}$. To keep the comparison causal, we choose the order on the 2020 period (25 weekly origins from 2 July to 17 December 2020) and evaluate it on the 2021 period (27 weekly origins from 15 April to 14 October 2021). The selection criterion is the mean origin-level 7-day MAPE of the linear forecaster on the smoothed target in 2020, minimized separately for each infectious period. As a check, we also use an aggregate criterion: the mean MAPE over both forecasters, both targets, and both horizons. We use infectious periods $D\in\{7,10,14\}$, horizons of 7 and 14 days, and two targets: the smoothed series $\nbar$ and raw reports $n$. Accuracy is measured by origin-level MAPE, MAE, and RMSE. Paired comparisons use the difference in origin-level MAPE, a HAC statistic with bandwidth one and a $t_{n-1}$ reference, where $n$ is the number of origins ($t_{24}$ in 2020 and $t_{26}$ in 2021) \cite{diebold1995,newey1987}, and Holm adjustment \cite{holm1979} within each family of 12 tests. As sensitivity checks we use bandwidths 2 and 3, and a noncircular moving-block bootstrap: overlapping blocks of 3 consecutive origins are drawn with replacement and concatenated to length $n$, with 10{,}000 resamples, 95\% percentile intervals for the mean difference, and random seed 90210.

\section{Results}\label{sec:results}
\subsection{Effect of the order on linear and constant-rate forecasts}
\begin{figure}[H]
\centering\includegraphics[width=0.95\linewidth]{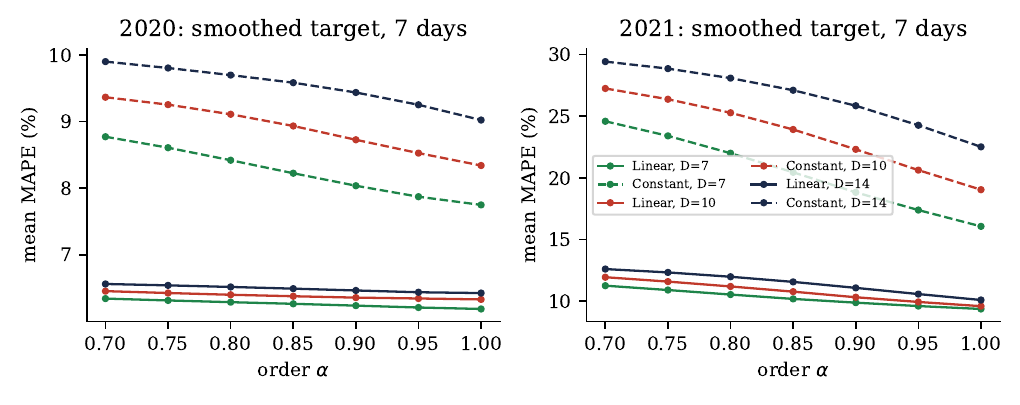}
\caption{Bangladesh: mean 7-day MAPE on the smoothed target against the fractional order, for the linear (solid) and constant-rate (dashed) forecasters and three infectious periods.}\label{fig:alpha}
\end{figure}

Table~\ref{tab:alpha} reports mean MAPE in 2021 for $D=14$ across all orders, and Figure~\ref{fig:alpha} shows the 7-day smoothed-target error for all infectious periods in both periods. In the 2021 results of Table~\ref{tab:alpha} the pattern is monotone. As $\alpha$ decreases from one, the error grows, for both forecasters, both targets, and both horizons. Across all 48 combinations of period, forecaster, target, horizon, and infectious period, the integer order gives the lowest mean MAPE in 47; the single exception (2020, linear, smoothed target, 14 days, $D=14$) favours $\alpha=0.95$ by 0.01 percentage points. Under MAE, $\alpha<1$ is best in 6 of 48 combinations, all in 2020 for the linear forecaster on the smoothed target, with gains of at most 1.6 cases per day. Under the selection criterion of Section~\ref{sec:design}, the order chosen on 2020 is $\alpha=1$ for every infectious period. The aggregate criterion makes the same choice (for $D=14$, a mean 2020 MAPE of 14.24\% at $\alpha=1$ against 14.92\% at $\alpha=0.70$). The single 2020 setting that favours $\alpha=0.95$ involves the 14-day horizon, which is not part of the selection criterion.

\begin{table}[H]
\centering\small
\caption{Mean MAPE (\%) over 27 origins in 2021 by fractional order, $D=14$.}\label{tab:alpha}
\begin{tabular}{lllrrrrrrr}
\toprule
Forecaster & Target & $H$ & 0.70 & 0.75 & 0.80 & 0.85 & 0.90 & 0.95 & 1.00\\
\midrule
Linear & Smoothed & 7 & 12.60 & 12.33 & 11.98 & 11.56 & 11.08 & 10.58 & \textbf{10.09}\\
Linear & Smoothed & 14 & 28.56 & 27.82 & 26.94 & 25.90 & 24.64 & 23.31 & \textbf{22.33}\\
Linear & Raw & 7 & 33.24 & 32.97 & 32.63 & 32.21 & 31.71 & 31.11 & \textbf{30.57}\\
Linear & Raw & 14 & 50.59 & 49.82 & 48.81 & 47.56 & 46.03 & 44.23 & \textbf{42.32}\\
Constant & Smoothed & 7 & 29.44 & 28.87 & 28.10 & 27.11 & 25.86 & 24.27 & \textbf{22.52}\\
Constant & Smoothed & 14 & 46.94 & 46.06 & 44.84 & 43.24 & 41.21 & 38.75 & \textbf{35.90}\\
Constant & Raw & 7 & 50.87 & 50.22 & 49.34 & 48.21 & 46.82 & 45.19 & \textbf{43.24}\\
Constant & Raw & 14 & 69.92 & 68.89 & 67.49 & 65.64 & 63.32 & 60.54 & \textbf{57.32}\\
\bottomrule
\end{tabular}
\end{table}

\subsection{Paired comparisons}
Table~\ref{tab:tests} compares $\alpha=0.9$ with $\alpha=1$ in 2021. For the linear forecaster, the fractional model is worse on average in all 12 comparisons, by 0.5 to 3.7 percentage points, but no difference is significant after Holm adjustment (smallest adjusted $p=0.20$). For the constant-rate forecaster, the fractional model is worse in all 12 comparisons, by 2.8 to 6.3 points, and all 12 differences remain significant after adjustment (adjusted $p\le0.011$). The HAC bandwidth does not change these conclusions: with bandwidths 2 and 3, again none of the linear and all 12 constant-rate comparisons are significant after Holm adjustment. The 95\% moving-block bootstrap intervals (unadjusted) exclude zero for all 12 constant-rate and 5 of 12 linear comparisons, always on the side favouring the integer order.

\begin{table}[H]
\centering\footnotesize
\caption{Difference in mean MAPE, $\alpha=0.9$ minus $\alpha=1$, over 27 origins in 2021 (percentage points; positive values favour the integer order), with Holm-adjusted HAC $p$-values in parentheses.}\label{tab:tests}
\begin{tabular}{llcccccc}
\toprule
 & & \multicolumn{2}{c}{$D=7$} & \multicolumn{2}{c}{$D=10$} & \multicolumn{2}{c}{$D=14$}\\
\cmidrule(lr){3-4}\cmidrule(lr){5-6}\cmidrule(lr){7-8}
Forecaster & Target & $H=7$ & $H=14$ & $H=7$ & $H=14$ & $H=7$ & $H=14$\\
\midrule
Linear & Smoothed & 0.51 (0.47) & 1.00 (0.73) & 0.72 (0.28) & 1.34 (0.73) & 0.98 (0.25) & 2.32 (0.40)\\
Linear & Raw & 0.75 (0.40) & 1.74 (0.73) & 0.82 (0.40) & 2.69 (0.30) & 1.14 (0.25) & 3.71 (0.20)\\
Constant & Smoothed & 2.76 (0.011) & 4.77 (0.011) & 3.28 (0.007) & 5.43 (0.011) & 3.34 (0.011) & 5.31 (0.011)\\
Constant & Raw & 3.64 (0.001) & 6.01 (0.007) & 3.92 (0.003) & 6.28 (0.009) & 3.58 (0.011) & 6.00 (0.011)\\
\bottomrule
\end{tabular}
\end{table}

\subsection{Sensitivity to the scaling of the removal rate}\label{sec:sens}
We repeated the linear and constant-rate analysis with $\gamma=(1/D)^\alpha$ instead of $\gamma=1/D$, for $\alpha=0.70$, $0.80$, $0.90$, $0.95$, and $1.00$. The two scalings coincide at $\alpha=1$. Under the alternative scaling, the 2020 period favours fractional orders: $\alpha<1$ has the lowest mean MAPE in 18 of 24 settings, and the selection criterion returns $\alpha=0.70$ for every infectious period. In 2021, however, $\alpha=1$ is best in all 24 settings, and the order selected on 2020 is worse than $\alpha=1$ in all 24 comparisons, for example 10.54\% against 9.37\% (linear, smoothed target, 7 days, $D=7$) and 23.81\% against 22.52\% (constant rate, smoothed target, 7 days, $D=14$). So the scaling changes which order fits the earlier period best, but under neither scaling did the order selected on 2020 outperform $\alpha=1$ in 2021.

\subsection{Neural forecaster}
For the neural ensemble we compare $\alpha=0.9$ with the integer order in both periods (Table~\ref{tab:neural}). Here the picture is different. The fractional model has the lower mean MAPE in 17 of 24 comparisons (9 of 12 in 2020 and 8 of 12 in 2021), for example 11.83\% against 13.42\% for the 7-day smoothed target with $D=14$ in 2021. The pattern is not consistent, however. With $D=10$ in 2021 the fractional model is worse in all four comparisons, by up to 6.7 percentage points. No difference is significant after Holm adjustment within either period (smallest adjusted $p=0.73$). All 1{,}560 neural fits at $\alpha=0.9$ converged without warnings. Which procedure is most accurate depends on the target and horizon. On the smoothed target the integer-order linear model remains the most accurate, for example 10.09\% against 11.83\% for the fractional neural model (7 days, $D=14$, 2021). On raw reports at 7 days, however, the neural forecasts of either order have lower mean MAPE than the linear model for every infectious period, for example 28.40\% (fractional) and 29.47\% (integer) against 30.57\% (linear) with $D=14$; in the integer-order case these neural-versus-linear differences are not significant after multiplicity adjustment \cite{hossain2026}.

\begin{table}[H]
\centering\small
\caption{Mean MAPE (\%) of the neural ensemble, $\alpha=0.9$ / $\alpha=1$; the lower value in each cell is in bold. No difference is significant after Holm adjustment.}\label{tab:neural}
\begin{tabular}{lllccc}
\toprule
Period & Target & $H$ & $D=7$ & $D=10$ & $D=14$\\
\midrule
2020 & Smoothed & 7 & \textbf{7.84} / 7.96 & \textbf{7.18} / 7.65 & \textbf{7.11} / 7.21\\
2020 & Smoothed & 14 & \textbf{12.98} / 15.75 & \textbf{12.08} / 13.12 & 12.39 / \textbf{11.86}\\
2020 & Raw & 7 & 15.91 / \textbf{15.72} & \textbf{16.03} / 16.14 & \textbf{15.95} / 16.02\\
2020 & Raw & 14 & \textbf{20.28} / 23.11 & \textbf{19.65} / 20.56 & 19.73 / \textbf{19.25}\\
\midrule
2021 & Smoothed & 7 & \textbf{11.60} / 12.85 & 12.25 / \textbf{10.86} & \textbf{11.83} / 13.42\\
2021 & Smoothed & 14 & \textbf{26.91} / 30.92 & 31.89 / \textbf{27.38} & \textbf{32.55} / 32.80\\
2021 & Raw & 7 & \textbf{28.94} / 29.66 & 28.79 / \textbf{26.83} & \textbf{28.40} / 29.47\\
2021 & Raw & 14 & \textbf{44.92} / 47.86 & 49.40 / \textbf{42.72} & \textbf{49.47} / 49.60\\
\bottomrule
\end{tabular}
\end{table}

\subsection{How memory changes the forecasts}
Why does memory raise the error of the linear and constant-rate forecasters? To look for an answer, Table~\ref{tab:bias} reports the mean signed percentage error of 7-day linear forecasts, split by the phase of the epidemic at the origin. The phase is defined by the model-based effective reproduction number $\mathcal{R}_{T-1}=\beta_{T-1}S_{T-1}/(\gamma N)$ of the integer-order reconstruction ($\alpha=1$, $\gamma=1/14$), so the same grouping of origins is used for both orders: growth if $\mathcal{R}_{T-1}>1$ (10 origins) and decline otherwise (17 origins). As a check that does not rely on the model, we also classify origins by the observed trend, with growth if $\nbar_T>\nbar_{T-7}$; the two classifications agree for 23 of 27 origins. In growing phases the two orders behave much the same. In declining phases, however, the fractional forecasts overshoot about twice as much, by 5.6\% against 2.9\%. The trend-based classification gives the same pattern (declining: 5.74\% against 2.87\%; growing: 1.89\% against 2.05\%). This pattern is consistent with a possible mechanism: through the terms $\sum_j|w_j|\hat I_{t+1-j}$ in \eqref{eq:prop}, the large prevalence of the preceding wave keeps feeding the forecast, so the fractional model may respond too slowly to a decline. The table describes an association and does not by itself establish this mechanism.

\begin{table}[H]
\centering\small
\caption{Mean signed percentage error (forecast minus target) of 7-day linear forecasts on the smoothed target, 2021, $D=14$, by phase at the origin, defined by $\mathcal{R}_{T-1}$ of the integer-order reconstruction.}\label{tab:bias}
\begin{tabular}{lrr}
\toprule
Phase & $\alpha=0.9$ & $\alpha=1$\\
\midrule
Growing ($\mathcal{R}_{T-1}>1$) & 2.18 & 1.97\\
Declining ($\mathcal{R}_{T-1}\le1$) & 5.57 & 2.92\\
\bottomrule
\end{tabular}
\end{table}

\subsection{Application to the United States}\label{sec:us}
We repeated the analysis for the United States with the same origins, training windows, forecasters, and settings, and widened the grid of orders to $\alpha\in\{0.30,0.35,\dots,1.00\}$. The reconstruction conserves the population to within $10^{-6}$ persons and keeps $\beta_t>0$ throughout, and all 4{,}680 linear and constant-rate forecast trajectories satisfy $S_t>0$ and $I_t>0$. Unlike in Bangladesh, susceptible depletion is substantial: at $\alpha=1$ the reconstructed $S_t/N$ falls to 0.826. For $\alpha<1$ it falls less, to 0.900 at $\alpha=0.9$, because the fractional recursion removes a fractional integral of the incidence rather than its sum.

\begin{figure}[H]
\centering\includegraphics[width=0.95\linewidth]{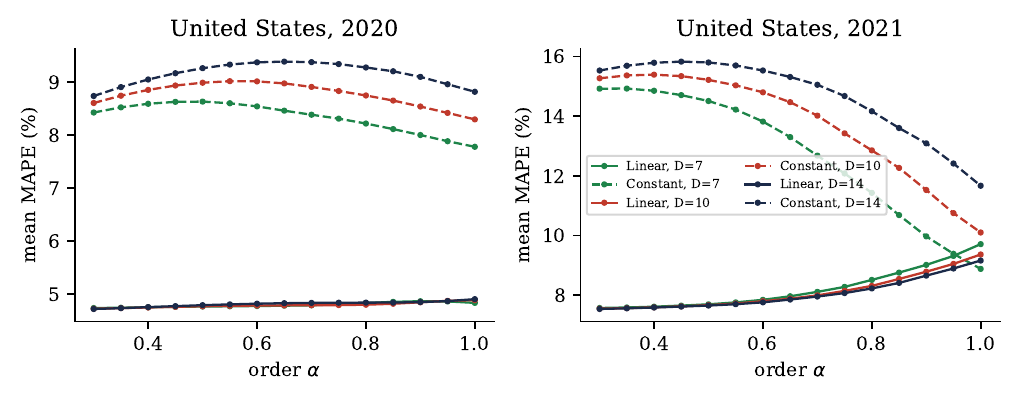}
\caption{United States: mean 7-day MAPE on the smoothed target against the fractional order, for the linear (solid) and constant-rate (dashed) forecasters and three infectious periods.}\label{fig:us}
\end{figure}

Table~\ref{tab:us} and Figure~\ref{fig:us} show a different pattern from Bangladesh. For the linear forecaster the error falls as $\alpha$ decreases, in every 2021 setting and on both targets. With $D=14$, for example, the 7-day smoothed-target MAPE is 7.54\% at $\alpha=0.30$ against 9.16\% at $\alpha=1$. The decrease continues down to the smallest order tested, so there is no interior optimum. The selection criterion of Section~\ref{sec:design} returns $\alpha=0.30$ for every infectious period; the aggregate criterion returns $\alpha=1$ for $D=7$ and $D=10$ and $\alpha=0.30$ for $D=14$. In 2021 the order chosen by the linear criterion lowers the linear forecaster's mean MAPE in all 12 comparisons, by 1.6 to 9.0 percentage points, but no difference is significant after Holm adjustment. For the constant-rate forecaster the same order raises the error in all 12 comparisons, by 1.5 to 8.5 points, again without significance after adjustment. Across all settings, the integer order gives the lowest mean MAPE in 22 of 48. The neural ensemble shows the reverse of the Bangladesh pattern: $\alpha=0.9$ gives lower mean MAPE than $\alpha=1$ in only 7 of 24 comparisons (3 of 12 in 2020 and 4 of 12 in 2021), and no difference is significant.

The bias analysis is consistent with damping. With $D=14$ in 2021, the integer-order linear forecasts overshoot by 4.47\% at origins in growing phases (13 origins). The overshoot falls to 2.37\% at $\alpha=0.9$ and becomes an undershoot of 3.26\% at $\alpha=0.30$, while errors in declining phases (14 origins) change little (1.73\%, 2.32\%, and 1.19\%). The constant-rate forecaster overshoots only slightly in growth (1.18\%) and gains nothing from the added inertia. Because daily US reports carry strong weekday artifacts, raw-target MAPE for the SIR-based forecasters ranged from about 65\% to 73\% with $D=14$ (about 63\% to 76\% across all infectious periods). For future seven-day case totals, which remove the weekday pattern, the linear forecaster's MAPE is 15.80\% at $\alpha=0.30$ and 16.07\% at $\alpha=1$, and the constant-rate forecaster's is 23.57\% and 17.32\%.

\begin{table}[H]
\centering\small
\caption{United States: mean MAPE (\%) over 27 origins in 2021 by fractional order, $D=14$; the lowest value in each row is in bold.}\label{tab:us}
\begin{tabular}{lllrrrrr}
\toprule
Forecaster & Target & $H$ & 0.30 & 0.50 & 0.70 & 0.90 & 1.00\\
\midrule
Linear & Smoothed & 7 & \textbf{7.54} & 7.65 & 7.96 & 8.66 & 9.16\\
Linear & Smoothed & 14 & \textbf{14.62} & 14.84 & 15.37 & 16.81 & 18.13\\
Linear & Raw & 7 & \textbf{65.27} & 65.41 & 65.88 & 66.63 & 67.49\\
Linear & Raw & 14 & \textbf{67.06} & 67.56 & 68.59 & 70.37 & 72.69\\
Constant & Smoothed & 7 & 15.53 & 15.80 & 15.06 & 13.09 & \textbf{11.67}\\
Constant & Smoothed & 14 & 23.96 & 24.12 & 23.06 & 20.83 & \textbf{19.25}\\
Constant & Raw & 7 & 66.38 & 66.34 & 65.73 & 65.08 & \textbf{64.91}\\
Constant & Raw & 14 & 71.85 & 71.79 & 70.97 & 70.23 & \textbf{69.74}\\
\bottomrule
\end{tabular}
\end{table}

\subsection{Seasonal naive benchmark and seven-day totals}\label{sec:totals}
The smoothed target averages seven daily reports, so at each origin its first six values contain reports that are already known when the forecast is made. We therefore also evaluate a target that lies entirely in the future: the total number of cases reported in the next seven days. We also add a seasonal naive benchmark that repeats the last seven reports observed at the origin. The forecast for day $T+h$ is the report for day $T+h-7$ when $h\le7$ and for day $T+h-14$ when $8\le h\le14$, so the last observed week is repeated twice. Persistence holds the last smoothed value fixed, so its seven-day total equals the total of the last seven reports, which is also the seasonal naive total. Table~\ref{tab:totals} reports the results for 2021 with $D=14$ for the linear and constant-rate forecasters.

On seven-day totals the linear forecaster is the most accurate method in both countries. The effect of memory is small and goes in opposite directions: in Bangladesh the error rises from 23.06\% at $\alpha=1$ to 24.41\% at $\alpha=0.9$, while in the United States it falls from 16.07\% to 15.80\% at $\alpha=0.30$. For the constant-rate forecaster, memory raises the error in both countries. On raw daily reports the two countries differ sharply. In Bangladesh the linear forecaster (30.57\% at seven days) beats seasonal naive (40.72\%). In the United States, seasonal naive reproduces the weekday pattern and reaches a seven-day raw MAPE of 28.33\%, against 64.91\% to 67.49\% for the SIR-based forecasters, which forecast a smooth level. When the target is daily reports, seasonal naive therefore performed better in the United States; the seven-day totals, evaluated separately above, give a different ranking.

\begin{table}[H]
\centering\small
\caption{Seven-day total forecast error (MAPE, \%) and raw daily MAPE of the seasonal naive benchmark, 2021, $D=14$; the lowest total-forecast error in each country is in bold.}\label{tab:totals}
\begin{tabular}{llrr}
\toprule
Method & Order & Bangladesh & United States\\
\midrule
\multicolumn{4}{l}{\emph{Seven-day totals}}\\
Linear & $\alpha=1$ & \textbf{23.06} & 16.07\\
Linear & $\alpha=0.9$ & 24.41 & 15.82\\
Linear & $\alpha=0.30$ & -- & \textbf{15.80}\\
Constant & $\alpha=1$ & 35.90 & 17.32\\
Constant & $\alpha=0.9$ & 39.76 & 19.47\\
Constant & $\alpha=0.30$ & -- & 23.57\\
Seasonal naive / persistence & -- & 34.56 & 18.91\\
\midrule
\multicolumn{4}{l}{\emph{Raw daily reports, seasonal naive}}\\
Seven days & -- & 40.72 & 28.33\\
Fourteen days & -- & 58.54 & 36.85\\
\bottomrule
\end{tabular}
\end{table}

\subsection{Comparison with standard baselines}\label{sec:baselines}
\begin{figure}[H]
\centering
\includegraphics[width=0.9\textwidth]{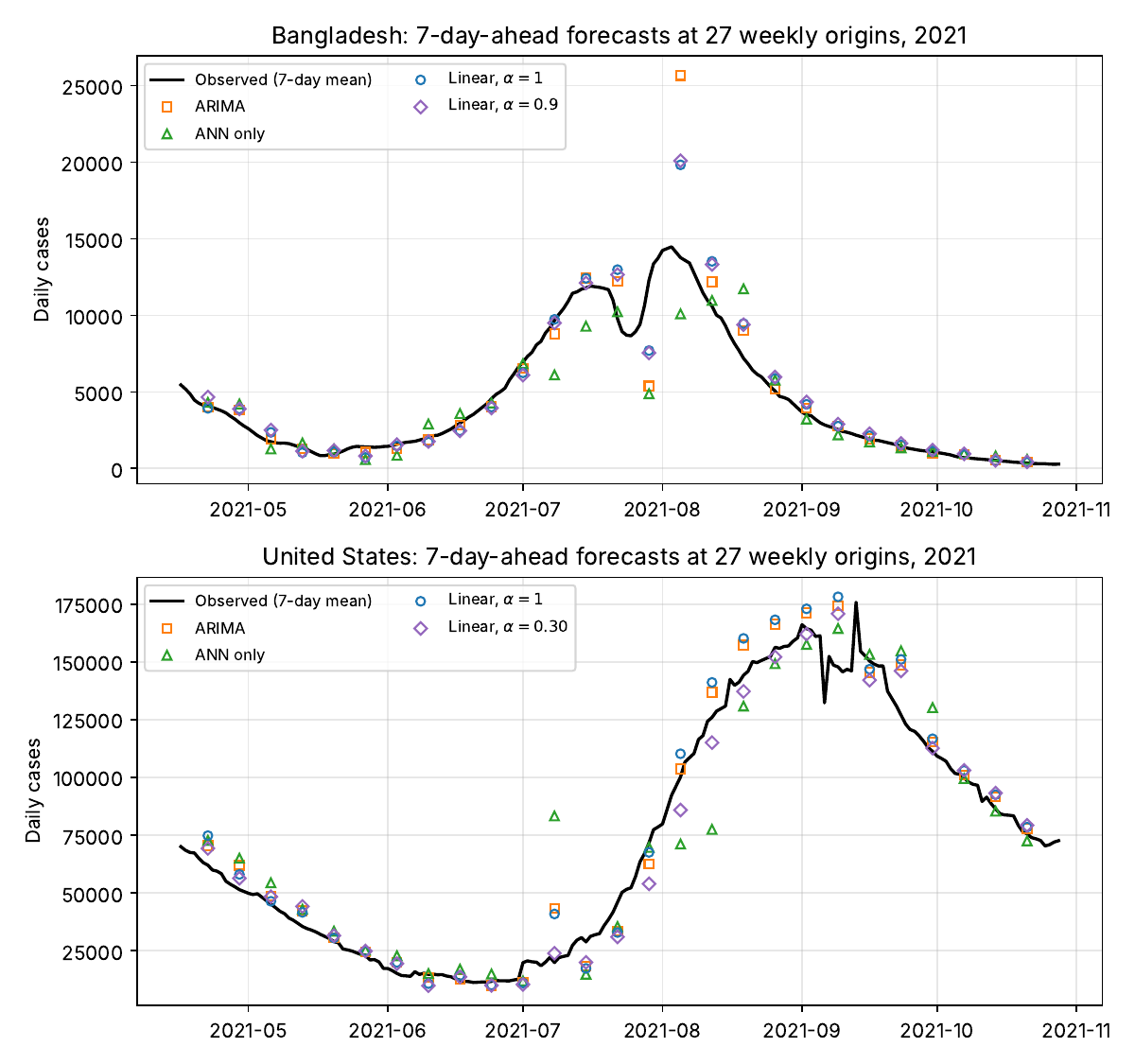}
\caption{Observed seven-day mean of daily cases (line) and seven-day-ahead forecasts made at the 27 weekly origins in 2021 (markers) by ARIMA, the ANN-only model and the linear SIR forecaster, $D=14$. Top: Bangladesh; bottom: United States.}\label{fig:baselines}
\end{figure}

Table~\ref{tab:baselines} compares the SIR-based forecasters with the baselines of Section~\ref{sec:baseline_models} on the smoothed target, and Figure~\ref{fig:baselines} shows the seven-day-ahead forecasts at every origin. Paired comparisons use the HAC test of Section~\ref{sec:design} with Holm adjustment within each family of three comparisons (main forecaster against each baseline, for a given country, target, horizon and metric). The main forecaster is the linear forecaster at the order selected on 2020: $\alpha=1$ in Bangladesh and $\alpha=0.30$ in the United States. ARIMA selected mostly low orders, ARIMA$(1,1,0)$ at 16 of the 27 Bangladesh origins and ARIMA$(1,1,1)$ at 14 of the 27 United States origins.

At the same order $\alpha=1$, forecasting the transmission rate improves on holding it fixed: the linear forecaster has a lower seven-day MAPE than the classical SIR model in Bangladesh (10.09\% against 22.52\%, Holm-adjusted $p=0.004$) and in the United States (9.16\% against 11.67\%, unadjusted HAC $p=0.084$). The selected United States forecaster ($\alpha=0.30$) reaches 7.54\% against 11.67\% for the classical SIR model ($p=0.026$); this difference reflects the overall performance of the two models, since both the transmission-rate forecast and the order change. The main linear forecaster ($\alpha=1$ in Bangladesh, $\alpha=0.30$ in the United States) is also more accurate on average than the ANN-only model (10.09\% against 14.54\% in Bangladesh, $p=0.063$; 7.54\% against 22.14\% in the United States, $p=0.46$), but this comparison changes both the forecasting method and the presence of SIR structure. The two neural forecasters use the same network architecture, which allows us to compare SIR-based transmission-rate forecasting with direct incidence forecasting. At seven days the neural SIR forecaster has the lower mean MAPE in both countries (13.42\% against 14.54\% in Bangladesh and 9.76\% against 22.14\% in the United States); at fourteen days it does so in the United States (17.47\% against 26.67\%) but not in Bangladesh (32.80\% against 31.56\%).

Against ARIMA the picture is mixed. In Bangladesh the linear forecaster has the lowest mean MAE and RMSE at both horizons, while ARIMA has the lowest mean MAPE (8.45\% against 10.09\% at seven days and 19.52\% against 22.33\% at fourteen days) and is the more accurate of the two at 18 of the 27 origins; the linear forecaster's advantage in MAE comes from two origins, 22 and 29 July 2021, when reports first fell sharply with a drop in testing around Eid al-Adha (21 July) and then rebounded; ARIMA extrapolated both the drop and the rebound more strongly than the SIR-based forecaster (Figure~\ref{fig:baselines}), and at the remaining 25 origins its MAE is lower on average. In the United States the fractional linear forecaster has the lowest mean MAE, RMSE and MAPE of all methods at both horizons (7.54\% against 9.14\% for ARIMA at seven days, and 14.62\% against 15.86\% at fourteen days). None of the differences from ARIMA is significant after Holm adjustment ($p\ge0.31$ in Bangladesh and $p\ge0.41$ in the United States). On raw daily reports (Table~\ref{tab:baselines_raw} in Appendix~B) all methods have much larger errors, and the ANN-only model is the least accurate in the United States.

\begin{table}[H]
\centering\small
\caption{Comparison with standard baselines: mean MAE and RMSE (cases per day) and MAPE over 27 weekly origins in 2021 on the smoothed target, $D=14$. The lowest value in each column and country is in bold.}\label{tab:baselines}
\resizebox{\textwidth}{!}{\begin{tabular}{lrrrrrr}
\toprule
 & \multicolumn{3}{c}{$H=7$} & \multicolumn{3}{c}{$H=14$}\\
\cmidrule(lr){2-4}\cmidrule(lr){5-7}
Model & MAE & RMSE & MAPE (\%) & MAE & RMSE & MAPE (\%)\\
\midrule
\multicolumn{7}{l}{\emph{Bangladesh}}\\
ARIMA & 515 & 637 & \textbf{8.45} & 1{,}223 & 1{,}488 & \textbf{19.52}\\
ANN only & 559 & 681 & 14.54 & 1{,}186 & 1{,}438 & 31.56\\
Constant rate, $\alpha=1$ (classical SIR) & 1{,}202 & 1{,}277 & 22.52 & 1{,}766 & 1{,}982 & 35.90\\
Linear, $\alpha=1$ & \textbf{445} & \textbf{558} & 10.09 & \textbf{1{,}097} & \textbf{1{,}351} & 22.33\\
Neural, $\alpha=1$ & 530 & 642 & 13.42 & 1{,}341 & 1{,}705 & 32.80\\
Linear, $\alpha=0.9$ & 464 & 581 & 11.08 & 1{,}134 & 1{,}396 & 24.64\\
Neural, $\alpha=0.9$ & 477 & 581 & 11.83 & 1{,}204 & 1{,}523 & 32.55\\
\midrule
\multicolumn{7}{l}{\emph{United States}}\\
ARIMA & 4{,}562 & 5{,}528 & 9.14 & 8{,}650 & 10{,}360 & 15.86\\
ANN only & 7{,}940 & 9{,}283 & 22.14 & 12{,}352 & 14{,}580 & 26.67\\
Constant rate, $\alpha=1$ (classical SIR) & 7{,}854 & 8{,}651 & 11.67 & 13{,}897 & 15{,}944 & 19.25\\
Linear, $\alpha=1$ & 4{,}926 & 5{,}928 & 9.16 & 10{,}143 & 12{,}272 & 18.13\\
Neural, $\alpha=1$ & 6{,}527 & 7{,}857 & 9.76 & 13{,}445 & 16{,}200 & 17.47\\
Linear, $\alpha=0.9$ & 4{,}715 & 5{,}679 & 8.66 & 9{,}567 & 11{,}535 & 16.81\\
Neural, $\alpha=0.9$ & 7{,}534 & 9{,}036 & 11.43 & 16{,}942 & 20{,}669 & 22.51\\
Linear, $\alpha=0.30$ & \textbf{4{,}039} & \textbf{4{,}963} & \textbf{7.54} & \textbf{7{,}822} & \textbf{9{,}512} & \textbf{14.62}\\
\bottomrule
\end{tabular}}
\end{table}

\section{Discussion}\label{sec:discussion}
The fractional model has memory in a precise sense. Theorem~\ref{thm:memory} shows that the unknown initial state is forgotten, but no faster than algebraically, and after 332 daily steps its influence is about $4.4\times10^{8}$ times larger at $\alpha=0.9$ than at $\alpha=1$; this is a property of the scheme and holds for any incidence series. Whether memory improves forecasts, however, depends on the setting. In Bangladesh, memory did not help the linear and constant-rate procedures: under the baseline scaling the integer order is best in 47 of 48 settings, a causal choice of $\alpha$ returns $\alpha=1$, and the loss from memory is large and significant for constant-rate forecasts. In the United States, memory lowered the errors of the linear forecaster in every 2021 comparison, although not significantly, while it raised those of the constant-rate forecaster.

The phase-specific errors are consistent with a damping effect, although changing $\alpha$ also changes the reconstructed series and its forecast propagation. In the United States the integer-order linear forecasts overshoot in growing phases, and memory reduces the overshoot. In Bangladesh the neural forecasts overshoot in growing phases \cite{hossain2026}, and memory lowered their error in 17 of 24 comparisons, whereas in the United States it did so in only 7 of 24. The absence of an interior optimum for the US linear forecaster supports this reading: there the order behaves like a damping parameter rather than an identified physical memory. None of the gains is statistically significant, so memory is not a reliable source of improvement in either country.

First, when $\beta$ and $\gamma$ are held fixed, the linearized threshold at disease-free states, $\beta S_*/(\gamma N)=1$, does not depend on the fractional order (Proposition~\ref{thm:stab}). What changes is how fast infectious perturbations decay below it. Because changing the order also changes the reconstruction and the training data, this slower decay may contribute to the larger errors we see during decline, but our experiments do not isolate it as the cause. Second, Theorem~\ref{thm:memory} matters beyond forecasting. Fractional epidemic models fitted to data keep the imprint of their initial conditions far longer than integer-order models, so the initial state should be estimated or varied in sensitivity analyses rather than fixed. Third, the reported-case series is already smoothed by a seven-day trailing average, which has a short memory of its own. This smoothing may partly explain why adding fractional memory did not consistently improve forecast accuracy, although we did not measure its effect separately.

The baseline comparison places these results in context. The reconstructed transmission rate is a useful quantity to forecast: forecasting it and propagating it through the SIR equations was more accurate on average than holding it constant, as in the classical SIR model. On the smoothed daily-case target, the neural SIR forecaster had lower mean MAPE than the ANN-only model in both countries at seven days and only in the United States at fourteen days. The SIR-based forecasters were competitive with ARIMA, the standard statistical benchmark, but did not significantly outperform it, and in Bangladesh ARIMA had the lower percentage error at most origins. With 27 origins per period the tests have limited power, so these comparisons establish competitiveness rather than superiority. The advantage of the SIR-based forecasters is that they produce epidemiologically interpretable quantities, such as the transmission rate and the susceptible and infectious compartments, which ARIMA does not.

The study has limitations. It uses two countries, treats reported cases as a proxy for infections, and for the United States relies on a national series with strong weekday artifacts; there, the best linear order lies at the edge of the grid, so the order is not identified. The removal parameter is fixed, and its reading as a mean infectious period holds only at $\alpha=1$. A single fractional order is used for all compartments. Theorem~\ref{thm:memory} bounds the rate of forgetting from below; it does not give the exact rate. The order was chosen from a grid rather than estimated together with the transmission rate, and other nonlocal operators, such as kernels with exponential or Mittag-Leffler tails, may behave differently. The Gr\"unwald--Letnikov scheme is first-order accurate, and higher-order methods such as the Adams predictor--corrector scheme \cite{diethelm2002} could confirm that the conclusions do not depend on the discretization. Future work could also produce calibrated prediction intervals and score them with proper scores such as the weighted interval score \cite{bracher2021}, estimate $\alpha$ together with an observation model, allow different orders for different compartments, and study diseases with longer infectious or reporting delays, such as dengue, where memory may matter more.

\section{Conclusion}\label{sec:conclusion}
We extended a hybrid SIR forecasting framework to Caputo fractional order. We proved well-posedness, conservation, and positivity, established a linearized threshold at disease-free states, and derived a discrete reconstruction that reduces exactly to the integer-order model. Our main theoretical result shows that for $\alpha<1$ initialization errors vanish, but no faster than algebraically, which gives memory a precise meaning.
In forecasting, the effect of memory depended on the country and the forecaster. In Bangladesh, under the baseline scaling of the removal rate, the integer-order model gave the lowest mean MAPE in 47 of the 48 settings involving linear and constant-rate forecasters, and selection using the 2020 data yielded $\alpha=1$; under the alternative scaling, the selected order ($\alpha=0.70$) had higher mean MAPE than $\alpha=1$ in all 24 comparisons in 2021. In the United States, a fractional order lowered the linear forecaster's error in all 12 comparisons in 2021 but raised the constant-rate forecaster's, and no difference was significant after Holm adjustment. The fractional neural forecaster had lower mean MAPE in 17 of 24 comparisons in Bangladesh and in 7 of 24 in the United States, again without significant differences. These results are consistent with memory damping over-responsive forecasts, although the source of this effect is not isolated. Against standard baselines, the SIR-based linear forecaster was more accurate than the classical SIR model and than a neural network without epidemic structure, and it was competitive with ARIMA without significantly outperforming it. The choice of a fractional order should therefore rest on out-of-sample evidence in each setting, not on in-sample fit or on the appeal of memory alone.

\section*{Appendix A. Sensitivity of the United States results to the 4 June 2021 correction}
Table~\ref{tab:usexcl} repeats the 2021 United States comparison between the order selected on 2020 ($\alpha=0.30$) and $\alpha=1$, leaving out every origin whose evaluation window contains the corrected value: the 3 June origin at seven days, and the 27 May and 3 June origins at fourteen days. Holm adjustment is applied within each forecaster's family of 12 tests.

\setcounter{table}{0}\renewcommand{\thetable}{A\arabic{table}}
\begin{table}[H]
\centering\footnotesize
\caption{United States, 2021, excluding origins affected by the correction: mean MAPE (\%) at $\alpha=0.30$ and $\alpha=1$, and the difference with its Holm-adjusted HAC $p$-value in parentheses; $n$ is the number of origins.}\label{tab:usexcl}
\begin{tabular}{lrrrrrrrrr}
\toprule
 & & & & \multicolumn{3}{c}{Linear} & \multicolumn{3}{c}{Constant}\\
\cmidrule(lr){5-7}\cmidrule(lr){8-10}
Target & $D$ & $H$ & $n$ & 0.30 & 1.00 & Diff. ($p$) & 0.30 & 1.00 & Diff. ($p$)\\
\midrule
Smoothed & 7 & 7 & 26 & 7.17 & 9.50 & -2.33 (1.00) & 14.87 & 8.81 & 6.06 (0.13)\\
Smoothed & 7 & 14 & 25 & 13.67 & 19.68 & -6.01 (1.00) & 22.67 & 14.77 & 7.90 (0.39)\\
Smoothed & 10 & 7 & 26 & 7.15 & 9.14 & -1.99 (1.00) & 15.20 & 10.11 & 5.09 (0.37)\\
Smoothed & 10 & 14 & 25 & 13.71 & 18.61 & -4.89 (1.00) & 23.07 & 17.09 & 5.98 (0.98)\\
Smoothed & 14 & 7 & 26 & 7.13 & 8.93 & -1.80 (1.00) & 15.45 & 11.77 & 3.68 (0.98)\\
Smoothed & 14 & 14 & 25 & 13.78 & 17.68 & -3.90 (1.00) & 23.36 & 19.65 & 3.71 (1.00)\\
Raw & 7 & 7 & 26 & 66.34 & 70.22 & -3.88 (1.00) & 66.94 & 65.12 & 1.82 (1.00)\\
Raw & 7 & 14 & 25 & 68.27 & 78.04 & -9.78 (1.00) & 71.65 & 69.26 & 2.40 (1.00)\\
Raw & 10 & 7 & 26 & 66.25 & 69.24 & -2.99 (1.00) & 67.09 & 65.46 & 1.63 (1.00)\\
Raw & 10 & 14 & 25 & 68.23 & 76.13 & -7.90 (1.00) & 71.86 & 70.25 & 1.61 (1.00)\\
Raw & 14 & 7 & 26 & 66.19 & 68.54 & -2.35 (1.00) & 67.21 & 65.93 & 1.28 (1.00)\\
Raw & 14 & 14 & 25 & 68.20 & 74.25 & -6.05 (1.00) & 72.03 & 71.34 & 0.69 (1.00)\\
\bottomrule
\end{tabular}
\end{table}

\section*{Appendix B. Baseline comparison on raw daily reports}
\setcounter{table}{0}\renewcommand{\thetable}{B\arabic{table}}

\begin{table}[H]
\centering\small
\caption{Comparison with standard baselines on raw daily reports: mean MAE and RMSE (cases per day) and MAPE over 27 weekly origins in 2021, $D=14$. The lowest value in each column and country is in bold.}\label{tab:baselines_raw}
\resizebox{\textwidth}{!}{\begin{tabular}{lrrrrrr}
\toprule
 & \multicolumn{3}{c}{$H=7$} & \multicolumn{3}{c}{$H=14$}\\
\cmidrule(lr){2-4}\cmidrule(lr){5-7}
Model & MAE & RMSE & MAPE (\%) & MAE & RMSE & MAPE (\%)\\
\midrule
\multicolumn{7}{l}{\emph{Bangladesh}}\\
ARIMA & 1{,}228 & 1{,}456 & 28.93 & 1{,}843 & 2{,}167 & \textbf{39.85}\\
ANN only & 1{,}189 & 1{,}398 & 31.81 & 1{,}660 & 1{,}905 & 47.95\\
Constant rate, $\alpha=1$ (classical SIR) & 1{,}711 & 1{,}877 & 43.24 & 2{,}181 & 2{,}468 & 57.32\\
Linear, $\alpha=1$ & 1{,}164 & 1{,}358 & 30.57 & 1{,}689 & 1{,}994 & 42.32\\
Neural, $\alpha=1$ & 1{,}059 & 1{,}227 & 29.47 & 1{,}738 & 2{,}081 & 49.60\\
Linear, $\alpha=0.9$ & 1{,}193 & 1{,}381 & 31.71 & 1{,}763 & 2{,}052 & 46.03\\
Neural, $\alpha=0.9$ & \textbf{987} & \textbf{1{,}179} & \textbf{28.40} & \textbf{1{,}593} & \textbf{1{,}904} & 49.47\\
\midrule
\multicolumn{7}{l}{\emph{United States}}\\
ARIMA & 27{,}000 & 33{,}479 & 69.69 & \textbf{28{,}831} & \textbf{35{,}936} & 72.03\\
ANN only & 31{,}031 & 36{,}914 & 103.47 & 33{,}306 & 39{,}985 & 93.01\\
Constant rate, $\alpha=1$ (classical SIR) & 27{,}611 & 34{,}118 & \textbf{64.91} & 30{,}939 & 37{,}893 & 69.74\\
Linear, $\alpha=1$ & \textbf{26{,}685} & \textbf{33{,}268} & 67.49 & 29{,}208 & 36{,}383 & 72.69\\
Neural, $\alpha=1$ & 27{,}656 & 33{,}816 & 66.01 & 31{,}685 & 38{,}824 & 69.15\\
Linear, $\alpha=0.9$ & 26{,}816 & 33{,}366 & 66.63 & 28{,}945 & 36{,}089 & 70.37\\
Neural, $\alpha=0.9$ & 28{,}215 & 34{,}039 & 65.77 & 33{,}925 & 40{,}920 & 73.17\\
Linear, $\alpha=0.30$ & 27{,}835 & 33{,}876 & 65.27 & 29{,}899 & 36{,}567 & \textbf{67.06}\\
\bottomrule
\end{tabular}}
\end{table}

\section*{Declarations}
\noindent\textbf{Data and code availability.} The compiled daily case record, the integer-order implementation of \cite{hossain2026}, whose independent reimplementation matched all 178 checked published values, and the complete fractional extension used here are provided in a reproducibility package. The package contains the scripts that produce every table and figure, the origin-level forecast metrics, the statistical-test and bootstrap outputs, the sensitivity analyses, the baseline comparison (ARIMA, ANN-only and classical SIR models, with the paired tests of Section~\ref{sec:baselines}), and the positivity checks of all 3{,}900 Bangladesh and 4{,}680 United States forecast trajectories; each check verifies $S_t>0$ and $I_t>0$ at every simulated step and records the minimum of $S_t/N$. The Bangladesh record is the compiled daily series of DGHS figures described in \cite{hossain2026}, and the United States series is from The New York Times \cite{nyt2021}. The package is available at \url{https://doi.org/10.5281/zenodo.24057230} (this DOI covers all versions of the package and resolves to the latest one).

\vspace{10pt}
\noindent\textbf{Funding}\\[3pt]
\noindent This research received no external funding.

\vspace{10pt}
\noindent\textbf{Competing interests}\\[3pt]
\noindent The authors declare no conflicts of interest regarding the publication of this paper.

\vspace{10pt}
\noindent\textbf{Use of generative AI}\\[3pt]
\noindent During the preparation of this manuscript, the authors used generative AI tools for code debugging and formatting. The authors reviewed all outputs and take full responsibility for the content.

\end{document}